\documentclass{amsart}
\usepackage{amsmath,amsthm}
\usepackage{graphics}
\usepackage{color}
\usepackage{epsfig}

\graphicspath{{./Figures/}}

\theoremstyle{plain}
\newtheorem{thm}{Theorem}[section]
\newtheorem{cor}[thm]{Corollary}

\newtheorem{prop}[thm]{Proposition}
\newtheorem{lemma}[thm]{Lemma}

\newtheorem*{thm4.1}{Theorem~\ref{thm:volumes}}
\newtheorem*{corblank}{Corollary}

\newcommand{\K}{\ensuremath{{\mathcal K}}}
\newcommand{\comment}[1]{}

\DeclareMathOperator{\vol}{vol}
\newcommand{\N}{\ensuremath{\mathbb{N}}}
\newcommand{\Q}{\ensuremath{\mathbb{Q}}}

\newcommand{\Z}{\ensuremath{\mathbb{Z}}}

\newcommand{\bdry}{\ensuremath{\partial}}

\begin{document}

%\title{Surgery descriptions and volumes of Berge knots.}
\title[Surgery descriptions and volumes of Berge knots I]{Surgery descriptions and volumes of Berge knots I:  Large volume Berge knots}

\author{Kenneth L. Baker}
\address{Department of Mathematics, University of Georgia \\ Athens, Georgia 30602}
\email{kb@math.uga.edu}

\thanks{This work was partially supported by a graduate traineeship from the VIGRE Award at the University of Texas at Austin and a VIGRE postdoc under NSF grant number DMS-0089927 to the University of Georgia at Athens.}



\begin{abstract}
By obtaining surgery descriptions of knots which lie on the genus one fiber of the trefoil or figure eight knot, we show that these include hyperbolic knots with arbitrarily large volume.  These knots admit lens space surgeries and form two families of Berge knots.  By way of tangle descriptions we also obtain surgery descriptions for these knots on minimally twisted chain links.
\end{abstract}

\maketitle

\section{Introduction}

In \cite{berge:skwsyls} Berge describes twelve families of knots that admit lens space surgeries.  These knots are referred to as {\em Berge knots} and appear to comprise all knots in $S^3$ known to have lens space surgeries.  Families (VII) and (VIII) are the knots which lie as essential simple closed curves on the fiber of the trefoil and figure eight knot respectively.  As perhaps suggested by the difference in the way Berge describes these two families as opposed to the other ten, these knots exhibit different properties with regards to hyperbolic volumes and accordingly surgery descriptions.

 In this article we show
\begin{thm4.1}
The knots which lie as simple closed curves on the fiber of the trefoil or figure eight knot (Berge's families (VII) and (VIII)) include hyperbolic knots with arbitrarily large volume.
\end{thm4.1}
By Thurston's Hyperbolic Dehn Surgery Theorem~\cite{thurston:gt3m}, this immediately implies the following corollary.
\begin{corblank}
There is no single link which admits surgery descriptions for every Berge knot.
\end{corblank}

In contrast, the sequel to this article \cite{baker:sdavobkII} shows that each Berge knot in the remaining families may be described by surgery on a minimally twisted five link of five components.  This implies the hyperbolic knots in these remaining families all have volume bounded above by the volume of a minimally twisted five chain link. 
 
There is, however, a unifying structure of surgery descriptions.   Theorem~\ref{thm:chainlinksurgdesc} shows that each knot in family (VII) or (VIII) may be described by surgery on a minimally twisted chain link of some number of components.  
%In contrast Theorem~\ref{thm:MT5Csurgdesc} shows that each knot in the remaining families may be described by surgery on a minimally twisted five link of five components.  As noted in Corollary~\ref{cor:bddvol}, this implies the hyperbolic knots in these remaining families all have volume bounded above by the volume of a minimally twisted five chain link.  Nevertheless, there is a unifying structure of surgery descriptions.   Theorem~\ref{thm:chainlinksurgdesc} shows that each knot in family (VII) or (VIII) may be described by surgery on a minimally twisted chain link of some number of components.  

\subsection{Acknowledgements}
The author wishes to thank John Luecke for his direction and many useful conversations.

\section{Definitions}

\subsection{Curves on tori and once punctured tori}
Let $\K$ be the set of isotopy classes of unoriented essential simple closed curves on the once-punctured torus $T$. 

Given two oriented simple closed curves $a$ and $b$ on $T$ that intersect once, we may orient $T$ so that $a \cdot b = +1$.  We say the ordered pair of oriented curves $\{a, b\}$ is a {\em basis} for $T$.  If $c$ is an essential (non-trivial, non-boundary parallel) simple closed curve on $T$ then, after choosing an orientation for $c$, $[c] = p [a] + q [b] \in H_1(T)$ for relatively prime integers $p$ and $q$.  Since $[-c] = -p[a]+-q[b]$, we may then identify the members $c \in \K$ with $\Q \cup \infty$ via the correspondence $c \mapsto p/q$.  The isotopy classes of unoriented non-trivial simple closed curves on the unpunctured torus are similarly identified with $\Q \cup \infty$ when given a basis.  These isotopy classes whether on the once-punctured or unpunctured torus are called {\em slopes}.  When given a basis, the corresponding rational number is also called the {\em slope} of the curve.

%%%%%
\subsection{Dehn surgery}
\begin{figure}[h]
\centering
\input{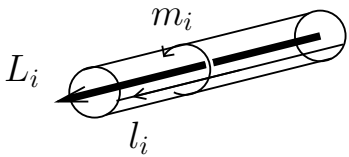}
\caption{The meridian and a longitude for $L_i$.}
\label{fig:std_basis}
\end{figure}

Given a link $L$ in an orientable $3$--manifold $M$, define the {\em link exterior} as $M_L = M - N(L)$.  For each component $L_i$ of $L$, let $\bdry_i M_L$ be the component of $\bdry M_L$ corresponding to $\bdry N(L_i)$.  Choose an orientation for $M$ and $L$.  Let $m_i$ be a simple closed curve on the torus $\bdry N(L)$ that bounds a disk in $\overline{N(L_i)}$ and is oriented such that the linking number of $m_i$ and $L_i$ is $+1$.  Let $l_i$ be an oriented simple closed curve on $\bdry L_i$ parallel to $L_i$.  Orient $\bdry N(L_i) = \bdry_i M_L$ so that $m_i \cdot l_i = +1$.  Note that this orientation on $\bdry_i M_L$ is induced from the orientation on $\overline{N(L_i)}$ rather than the orientation on $M_L$.  
The slope $m_i$ is called the {\em meridian} of $L_i$, and the slope $l_i$ is called a {\em longitude} of $L_i$.  Taken together, $\{m_i, l_i\}$ forms a basis for $\bdry_i M_L$.  Figure~\ref{fig:std_basis} shows $L_i$, $m_i$, and $l_i$.

Assume $L$ is a $k$ component link.  Let $\rho_i$ be a slope on $\bdry_i M_L$, and let $\overline{\rho} = (\rho_1, \dots, \rho_k)$.  Then the {\em $\overline{\rho}$--Dehn surgery} on $L$ in $M$ or  the {\em $\overline{\rho}$--Dehn filling} of $M_L$ is the glueing of a solid torus onto  each $\bdry_i M_L$ such that $\rho_i$ bounds a disk in the attached solid torus.  We may choose to set $\rho_i = *$ if the component $L_i$ or $\bdry_i M_L$ is to be left unsurgered or unfilled respectively.

\subsection{Dehn twists}

\begin{figure}[h]
\centering
\input{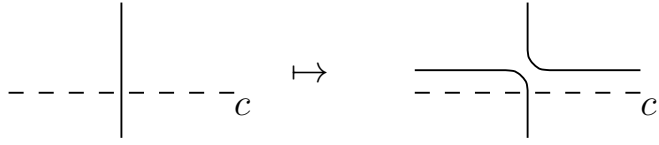}
\caption{Left-handed Dehn twist about $c$.}
\label{fig:Dehn_twist}
\end{figure}

Let $c$ be a simple closed curve on an oriented surface $S$.  
The homeomorphism $S \to S$ defined by cutting $S$ along $c$ and regluing along $c$ after rotating a full counter-clockwise turn is the {\em left-handed Dehn twist} of $S$ about $c$. 
Figure~\ref{fig:Dehn_twist} illustrates this homeomorphism.  This homeomorphism is the identity outside a neighborhood of $c$.
The handedness of the Dehn twist is observed from the positive side of $S$.  A left-handed Dehn twist when viewed from the negative side of $S$ appears  right-handed.

Given an annulus $R$ embedded in an orientable $3$-manifold $M$, set $\bdry R = L_{+1} \cup L_{-1}$.  Orient $L_{+1}$ and $L_{-1}$ coherently.  Form the bases $\{m_i, l_i\}$, $i=\pm1$, where $m_i$ is the standard meridian on $\bdry N(L_i)$ and $l_i$ is the longitude on $\bdry N(L_i)$ induced by $R$ for $i = \pm 1$.  (That is, let $l_i = \bdry N(L_i) \cap R$.)

\begin{figure}[h]
\centering
\input{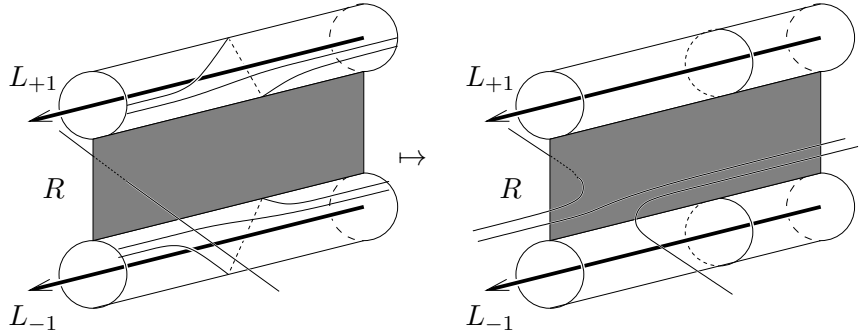}
\caption{Twisting along an annulus $R$ by surgery.}
\label{fig:surgery_twist}
\end{figure}

In these bases, Dehn surgery of $\frac{1}{n}$ on $L_{+1}$ and $-\frac{1}{n}$ on $L_{-1}$ yields the homeomorphism $M \cong M_{(L_{+1} \cup L_{-1})}(\frac{1}{n},-\frac{1}{n})$.  This homeomorphism is realized by cutting $M_{(L_{+1} \cup L_{-1})}$ along $R$, then regluing along $R$ after spinning $n$ times counter-clockwise, and filling the resulting two toroidal boundary components incident to $R$ trivially.  See Figure~\ref{fig:surgery_twist} for the case $n=-2$ and its effect upon a curve that intersects $R$.  Note that this homeomorphism is the identity outside a neighborhood of $R$.

If an orientable surface $S$ intersects the annulus $R$ along a simple closed curve $c$ which is essential in $R$, then the homeomorphism \[M \to M_{(L_{+1} \cup L_{-1})}(\frac{1}{n},-\frac{1}{n})\] restricted to $S$ is the composition of $n$ Dehn twists.  If $n > 0$ and we are viewing $S$ from the same side as $L_{+1}$, these are left-handed Dehn twists.

%%%%%%%%%
\subsection{Continued fractions and curves on once-punctured tori.}
We define a {\em continued fraction expansion} $[r_n, r_{n-1}, \dots, r_2, r_1]$ of $p/q \in \Q \cup \infty$ of {\em length} $n \in \N$ as follows:
\[ \frac{p}{q} = [r_n, r_{n-1}, \dots, r_2, r_1] = 
\cfrac{1}{r_n-
 \cfrac{1}{r_{n-1}-
  \cfrac{1}{\dots -
   \cfrac{1}{r_2 - 
    \cfrac{1}{r_1
}}}}}
\]
where $r_i \in \Z$.  
The continued fraction of length $0$ is $[\emptyset]$ and represents $\infty$.  Note that a rational number has many continued fraction expansions.  
%%%%

Fix a basis $\{a, b\}$ for the once punctured torus $T$.  
Let $\alpha$ and $\beta$ be the left-handed Dehn twists of $T$ along $a$ and $b$ respectively.
 
%%%
\begin{lemma}\label{lem:KandCFE} %different from the one in diss
A curve $K \in \K$ has slope $[r_n, \dots, r_2, r_1]$, $n$ odd, if and only if
\[
K = \beta^{r_n} \circ \cdots \circ \alpha^{r_2} \circ \beta^{r_1} (a).
\]
\end{lemma}

\begin{proof}
Recall if $p/q$ is the slope of the curve $K \in \K$, then, for some orientation on $K$, $[K] = p [a] + q[b] \in H_1(T)$ assuming $p$ and $q$ coprime.  One may check that $[\alpha^r (K)] = (p-r q)[a] + q[b]$ yields the slope $(p-r q)/q$ for $\alpha^r (K)$ and that $[\beta^r (K)] = p[a] + (r p+q)[b]$ yields the slope $p/(r p+q)$ for $\beta(K)$.
Performing these Dehn twists in succession, 
\[ [\beta^t \circ \alpha^s (K)] = (p-sq)[a] + (t(p-sq)+q)[b] = (p-sq)[a] + (tp + (1-ts)q)[b] \]
so that $\beta^t \circ \alpha^s (K)$ has slope 
\[\frac{p-sq}{t(p-sq)+q} = 
\cfrac{1}{t + \frac{q}{p-sq}} = \cfrac{1}{t - \cfrac{1}{s - \frac{p}{q}}}.\]

Given a curve $K \in \K$ with slope $[r_n, \dots, r_2, r_1]$, $n$ odd, let $K' \in \K$ be the curve with slope $[r_{n-2}, \dots, r_2, r_1]$.  Since
\[[r_n, \dots, r_2, r_1] = \cfrac{1}{r_n - \cfrac{1}{r_{n-1} - [r_{n-2}, \dots, r_2, r_1]}},\]
one concludes that $K = \beta^{r_n} \circ \alpha^{r_{n-1}} (K')$.  Similarly, given that $K = \beta^{r_n} \circ \alpha^{r_{n-1}} (K')$ and that $K'$ has slope $[r_{n-2}, \dots, r_2, r_1]$, one concludes that $K$ has slope $[r_n, \dots, r_2, r_1]$.
With these observations and that $[r]$ is the slope of $\beta^r(a)$, the lemma follows from an induction argument.
\end{proof}

Note 
\begin{enumerate}
\item $[r_n, \dots, r_2, r_1] = [r_n, \dots, r_2, r_1 \pm 1, \pm 1]$ so that any continued fraction of even length may be altered to a continued fraction of odd length, and
\item if 
$K = \alpha^{r_n} \circ \cdots \circ \alpha^{r_2} \circ \beta^{r_1} (a)$,
then 
$K = \beta^0 \circ \alpha^{r_n} \circ \cdots \circ \alpha^{r_2} \circ \beta^{r_1} (a)$
so that slope of $K$ has the continued fraction expansion $[0, r_n, \dots, r_2, r_1]$ of odd length.
\end{enumerate}
%%%%%%%

\subsection{Tangles}
 A {\em tangle} $(B,t)$ is a pair consisting of a punctured $3$-sphere $B$ and a properly embedded collection $t$ of disjoint arcs and simple closed curves.  Two tangles $(B_1, t_1)$ and $(B_2, t_2)$ are homeomorphic if there is a homeomorphism of pairs
\[h \colon (B_1, t_1) \to (B_2, t_2). \]

A boundary component $(\bdry B_0, t \cap \bdry B_0)$ of a tangle $(B,t)$ is a sphere $\bdry B_0$ together with some finite number of points $p=t \cap \bdry B_0$.  Typically we consider the situation where $p$ consists of just four points.

Given a sphere $S$ and set $p$ of four distinct points on $S$,
 a {\em framing} of $(S, p)$ is an ordered pair of (unoriented) simple closed curves $(\hat{m}, \hat{l})$ on $S - N(p)$ such that each curve separates different pairs of points of $p$.
 
\begin{figure}
\centering
\input{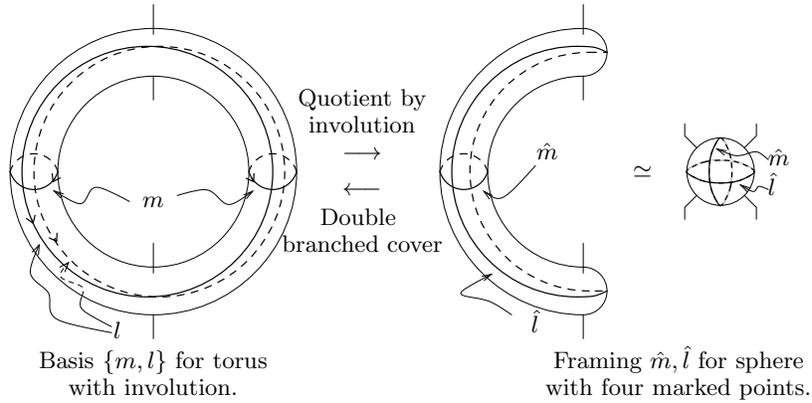}
\caption{The correspondence of bases and framings.}
\label{fig:tangleframing}
\end{figure}

 Let $(S, p)$ be a sphere with four points with framing $(\hat{m}, \hat{l})$.  The double cover of $S$ branched over $p$ is a torus.  Single components, say $m$ and $l$, of the lifts of the framing curves $\hat{m}$ and $\hat{l}$ when oriented so that $m \cdot l = +1$ (with respect to the orientation of the torus) form a basis for the torus.  Similarly, a basis on a torus induces a framing on the sphere with four points obtained by the quotient of an involution that fixes four points on the torus.  See Figure~\ref{fig:tangleframing}.
 
\section{Surgery descriptions of knots on a once-punctured torus fiber}

\begin{figure}[h]
\begin{center}
\input{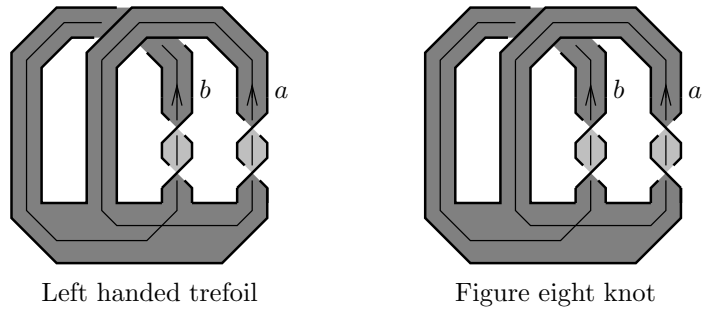}
\end{center}
\caption{Left handed trefoil and figure eight knot each with fiber and basis.}
\label{fig:genonefiberedknots}
\end{figure}

We develop a surgery description of both knots on the fiber of the trefoil and knots on the fiber of the figure eight knot.  A fiber together with a basis are shown in Figure~\ref{fig:genonefiberedknots} for each the trefoil and figure eight knot.

\comment{
\begin{figure}[h]
\begin{center}
\input{Figures/LH_Trefoil_fiber.pstex_t}
\end{center}
\caption[Left handed trefoil]{Left handed trefoil with fiber and basis.}
\label{fig:+tref}
\end{figure}

\begin{figure}[h]
\centering
\input{Figures/FE_Knot_fiber.pstex_t}
\caption[Figure eight knot]{Figure eight knot with fiber and basis.}
\label{fig:fig8}
\end{figure}
}

Let $M_{\eta}$ denote the fiber bundle $T \times [-n,n+1] / \eta$ where ${(x, n+1)\sim(\eta(x),-n)}$ for some $n \in \N$ and some homeomorphism $\eta \colon T \to T$.  It will be convenient to work with $\eta^{-1}$ instead of $\eta$.  If $\eta^{-1} = \beta \circ \alpha$, then $M_{\eta}$ may be recognized as the complement of the left handed trefoil.  If $\eta^{-1} = \beta \circ \alpha^{-1}$, then $M_{\eta}$ is the complement of the figure eight knot. 
Note that in either case $M_\eta \subset S^3$.
Throughout, $\eta^{-1}$ is assumed to be either $\beta \circ \alpha$ or $\beta \circ \alpha^{-1}$.

Consider the collection $\K$ of essential simple closed curves on the once-punctured torus $T$.  Let $\K_\eta$ denote the image of $\K$ under the inclusion
\[T \to T \times \{0\} \subset T \times I /_\eta = M_\eta \subset S^3.\]
Thus depending on $\eta$, $\K_\eta$ is assumed to be either the collection of knots that lie as essential simple closed curves on the fiber of the left handed trefoil or those that lie on the fiber of the figure eight knot.  Up to mirror image, these form Berge's families VII and VIII and are the primary subjects of our study.

\subsection{The link $L(2n+1, \eta)$ and the knots of $\K_\eta$}

Given $M_\eta = T \times [-n, n+1]/_\eta \subset S^3$ with $n \in \N$, set $L_i = a \times\{i\}$ or $L_i = b \times \{i\}$ if $i$ is even or odd respectively, $i \in \{-n, -n+1, \dots , n\}$.  Define $L(2n+1, \eta)$ to be the oriented link 
\[L(2n+1, \eta) = \bigcup_{i=-n}^{n} L_i \subset M_{\eta} \subset S^3. \]
We write simply $L(2n+1)$ when $\eta$ is understood.

\begin{figure}
\centering
\input{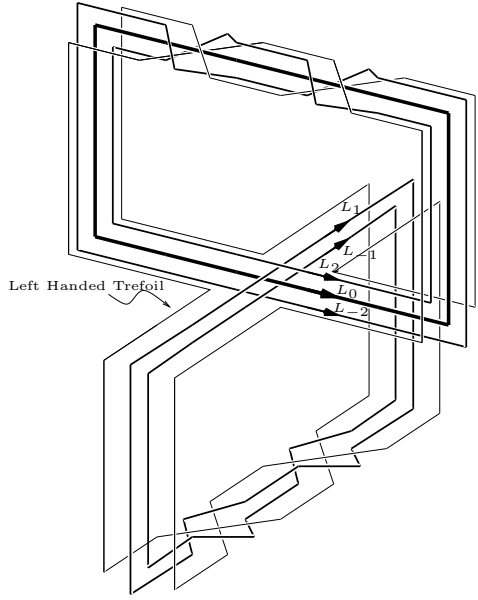}
\caption{The link $L(5)$ and the left handed trefoil.}
\label{fig:LinkLHT}
\end{figure}

Figure~\ref{fig:LinkLHT} shows the link $L(5)$ where $\eta^{-1} = \beta \circ \alpha$.  Included in the figure for reference is the corresponding genus one fibered knot, the left handed trefoil.

%\subsection{Obtaining $K \in \K_\eta$ by surgery on $L(2n+1,\eta)$}

Consider the exterior of $L(2n+1)$, $S^3_{L(2n+1)}$.  For each $i \in \{-n, \dots, n\}$, let $\mu_i$ be the standard meridian for $\bdry N(L_i)$, and let $\lambda_i$ be the longitude associated to the isotopy class of a component of $\bdry N(L_i) \cap T \times \{i\}$. 

Notice that $\lambda_i$ is not the standard longitude for $L_i$ as the two have non-zero linking number.  Moreover, for each integer $0 < i \leq n$, the representatives of $\lambda_i$ and $\lambda_{-i}$ together bound an annulus $R_i$ in $S^3_{L_i \cup L_{-i}}$, the exterior of $L_i \cup L_{-i}$, that intersects $T \times \{0\}$ along a curve isotopic to either $b \times \{0\}$ if $i$ is odd or $a \times \{0\}$ if $i$ is even.  Therefore, with respect to the bases $\{\mu_i, \lambda_i\}$ and $\{\mu_{-i}, \lambda_{-i}\}$, the surgery $S^3_{L_{-i} \cup L_i} (-\frac{1}{r_i}, \frac{1}{r_i})$ realizes a ``spin'' along $R_i$.  On $T \times \{0\}$ this restricts to the Dehn twist $\alpha^{r_i}$ or $\beta^{r_i}$ depending on the parity of $i$.  

\begin{prop}\label{prop:K_as_surgery}
If $K \in \K_\eta$, then $K$ may be obtained by surgery on ${L(2n+1,\eta)}$ for some $n \in \N$.
\end{prop}

\begin{proof}
With respect to the basis $\{a, b\}$ on $T$, $K$ has the slope $p/q$ which in turn has a continued fraction expansion $[r_n, \dots, r_2, r_1]$ of odd length.
By Lemma~\ref{lem:KandCFE}, this continued fraction expansion implies an expression for the knot $K \subset T \times \{0\}$ as the image of a composition of Dehn twists:
\[K = \beta^{r_n} \circ \cdots \circ \alpha^{r_2} \circ \beta^{r_1} (a).\]

By nesting the surgery realizations of the Dehn twists, we may obtain $K$ by surgery on $L(2n+1)$.
Let $\rho_i =\frac{1}{r_i}$ or $-\frac{1}{r_{-i}}$ if $i$ is positive or negative respectively.
Set $\rho_0 = *$ to indicate that Dehn surgery is not done on $L_0$.  Let $\overline{\rho} = (\rho_{-n}, \rho_{-n+1}, \dots, \rho_{n-1}, \rho_n)$.  Then $K$ is the image of $L_0$ in $S^3_{L(2n+1)}(\overline{\rho})$.
\end{proof}

\section{Volumes}

\begin{thm}\label{thm:volumes}
The knots which lie as simple closed curves on the fiber of the trefoil or figure eight knot (Berge's families (VII) and (VIII)) include hyperbolic knots with arbitrarily large volume.
\end{thm}

In the proof of  Theorem~\ref{thm:volumes}, we employ Corollary~\ref{cor:Lishyp} which states that for $n>1$, $L(2n+1,\eta) \subset S^3$ is a hyperbolic link.  We defer this corollary and its proof until Section~\ref{sec:HoLE}.

\begin{proof}
This theorem concerns the knots in $\K_\eta$ for $\eta$ corresponding to each the trefoil an the figure eight knot.
By Proposition~\ref{prop:K_as_surgery}, we may obtain these knots as surgeries on $L(2n+1, \eta)$ according to their description as a product of Dehn twists.  

By Corollary~\ref{cor:Lishyp}, $S^3_{L(2n+1)}$ is a hyperbolic manifold with $2n+1$ cusps for $n \geq 2$.  Adams~\cite{adams:volumes} shows that a hyperbolic manifold with $N$ cusps has volume at least $v_N \geq n \cdot v$.  Here $v_N$ is the volume of the smallest volume $N$ cusped hyperbolic manifold, and $v$ is the volume of a regular ideal tetrahedron.  Thus $\vol(S^3_{L(2n+1)}) \geq (2n+1) \cdot v$.  Thurston's Hyperbolic Dehn Surgery Theorem~\cite{thurston:gt3m} states that large enough filling on a cusp of a (finite volume) hyperbolic manifold will yield a hyperbolic manifold with volume near that of the unfilled manifold.  
Thus given $\varepsilon > 0$ we may choose $|r_1|, \dots, |r_n| \gg 0$ so that if $K \in \K_\eta$ corresponds to the continued fraction expansion $[r_1, r_2, \dots, r_k]$ as in Lemma~\ref{lem:KandCFE} and $\overline{\rho}$ is chosen as in Proposition~\ref{prop:K_as_surgery}, then 
\[\vol(K) = \vol(X_{L(2n+1)}(\overline{\rho}) - N(L_0)) > \vol(X_{L(2n+1)}) - \varepsilon.\]
Hence $\vol(K) > (2n+1) \cdot v - \varepsilon$.

For each $n > 1$ we may thus choose as sequence of knots $\{K_n\}_{n=2}^{\infty} \subset \K_\eta$ with corresponding continued fraction expansions $[r_1(n), \dots, r_n(n)]$ where for each $n$ the $|r_i(n)|$ are sufficiently large so that $\vol(K_n) \to \infty$ as $n \to \infty$.
\end{proof}

\section{Homeomorphism of Link Exteriors}\label{sec:HoLE}

\begin{figure}[h!]
\begin{center}
\input{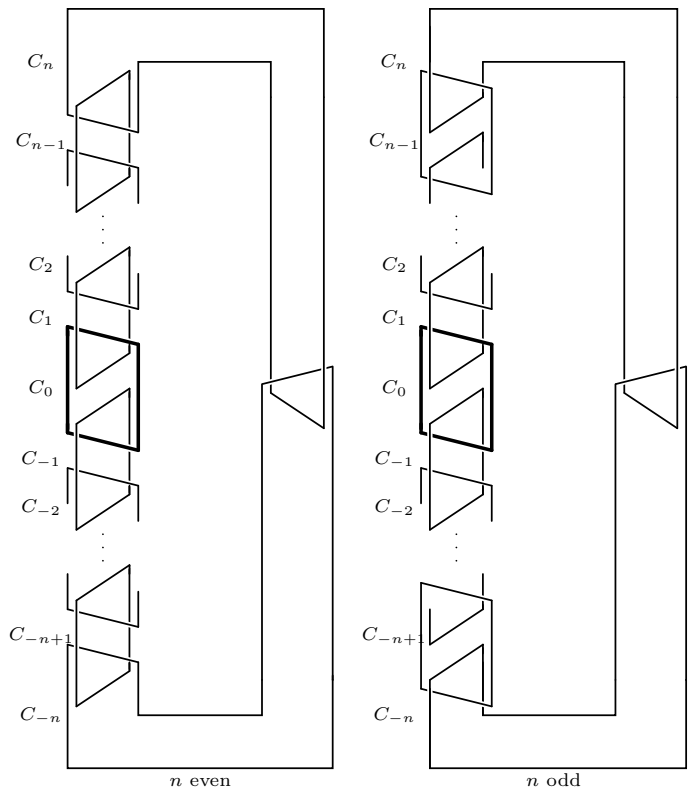}
\end{center}
\caption[Minimally twisted $2n+1$ chain.]{The minimally twisted $2n+1$ chain for $n$ even and $n$ odd.}
\label{fig:MT2n+1C}
\end{figure}

For $n \in \N$, let $C(2n+1) = \overset{n}{\underset{i=1}{\cup}} C_i$ be the {\em minimally twisted chain link of $2n+1$ components} as shown for each $n$ even and $n$ odd in Figure~\ref{fig:MT2n+1C}.  Note there are actually two such links, $C(2n+1)$ and its reflection $-C(2n+1)$, depending on how the clasping of $C_{n}$ and $C_{-n}$ is done.  For each $i \in \{-n, \dots, n\}$, let $\{m_i, l_i\}$ be the standard meridian-longitude pair for $C_i$.

The main theorem of this section is
\begin{thm}\label{thm:homeomorphism}
The links $L(2n+1,\eta)$ and $C(2n+1)$ have homeomorphic exteriors.  I.e.\ $S^3_{L(2n+1,\eta)} \cong S^3_{C(2n+1)}$.
\end{thm}
Consequentially, we will use this homeomorphism to obtain descriptions of the knots $K \in \K_\eta$ as surgeries on these chain links.  Through different methods Yuichi Yamada~\cite{yy} has also obtained descriptions of many of these knots as surgeries on minimally twisted chain links.  In the following section we describe the other families of Berge knots as surgeries on the minimally twisted five chain link $C(5)$.
Furthermore, though it may be done directly, this homeomorphism simplifies the proof of the hyperbolicity of $L(2n+1)$.

\subsection{Proof of Theorem \ref{thm:homeomorphism}} 

\begin{figure}[h]
\centering
\input{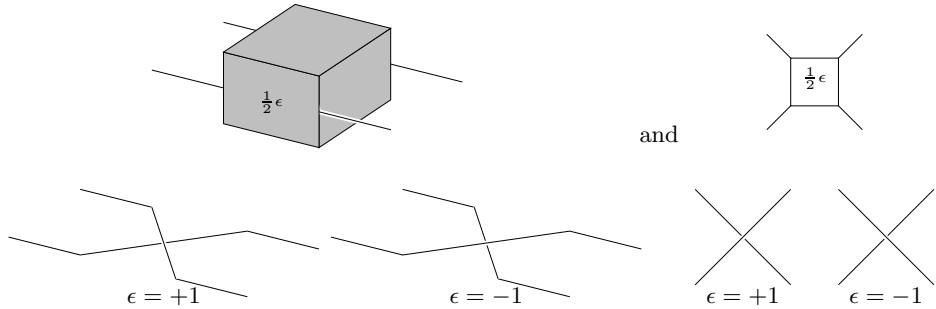}
\caption{Twistings of $\frac{1}{2} \epsilon$ for $\epsilon = \pm1$.}
\label{fig:epsilonlegend}
\end{figure}

Throughout this subsection and the next, the choice of $\epsilon = \pm 1$ depends on the choice of $\eta$.  This is encapsulated by 
\[ \eta^{-1} = \beta \circ \alpha^{\epsilon} \]
so that $\epsilon = +1$ corresponds to the left hand trefoil and $\epsilon = -1$ corresponds to the figure eight knot.  Furthermore, in the ensuing figures, blocks with $\frac{1}{2} \epsilon$ indicate the half twists as shown in Figure~\ref{fig:epsilonlegend}.

\begin{proof}
The links $L(2n+1,\eta)$ all admit a strong involution as shown in Figure~\ref{fig:linkinvolution} for $n=2$.  The genus one fibered knot is included for reference as it is invariant under the involution.  Quotienting the link exterior $S^3_{L(2n+1,\eta)}$ by the involution yields the tangle $(B_{2n+1},t_L)$ where $B_{2n+1}$ is a $2n+1$ punctured $S^3$ and $t_L$ is the image of the involution axis.  Hence the double branched cover of $B_{2n+1}$ branched over $t_L$ is $S^3_{L(2n+1, \eta)}$.  A sequence of isotopies of $(B_{2n+1}, t_L)$ for $n = 2$ into a ``nice'' form are indicated in Figure~\ref{fig:linkquotient}.  The quotient of the genus one fibered knot is also shown in these figures.

\begin{figure}
\centering
\input{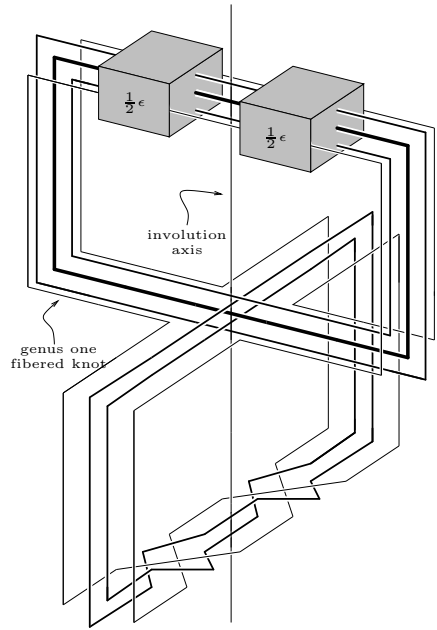}
\caption{The link $L(5, \eta)$ with an axis of strong involution and the genus one fibered knot.}
\label{fig:linkinvolution}
\end{figure}

\begin{figure}
\centering
\input{Figures/GenOne-isotopycompl-comb.pstex_t}
\caption{The quotient of the complement of the link $L(5, \eta)$ by the strong involution, $(B_5,t_L)$, followed by isotopies.}
\label{fig:linkquotient}
\end{figure}

\comment{
\begin{figure}
\centering
\input{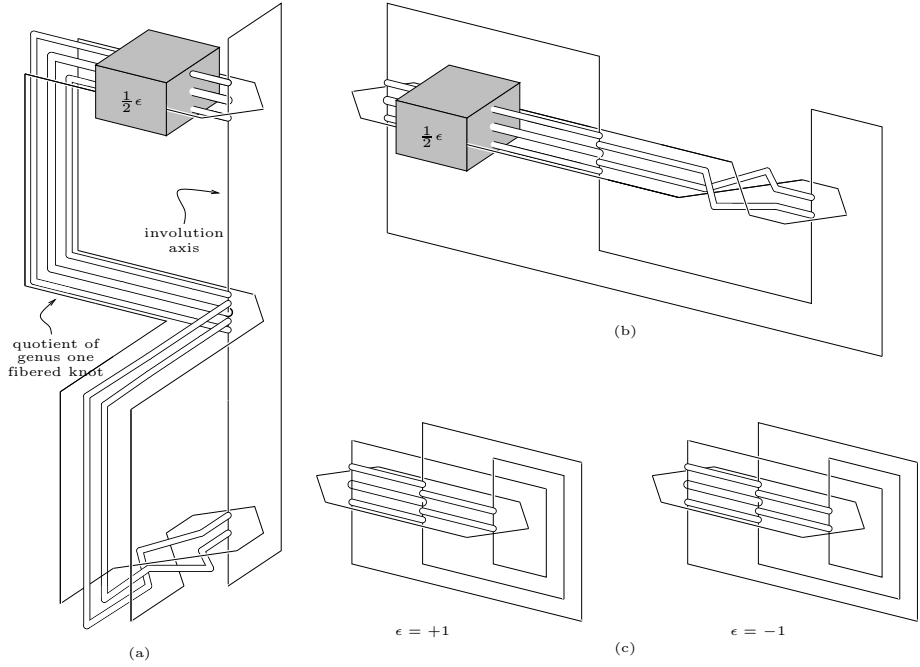}
\caption{The quotient of the complement of the link $L(5, \eta)$ by the strong involution, $(B_5,t_L)$, followed by isotopies.}
\label{fig:linkquotient}
\end{figure}

\begin{figure}
\centering
\input{Figures/GenOne-isotopycompl2.pstex_t}
\caption{Further isotopies of $(B_5,t_L)$.}
\label{fig:linkquotient2}
\end{figure}
}

We may trace the longitudes $\lambda_i$ and standard meridians $\mu_i$ of $\bdry N(L_i)$ through the quotient of $S^3_{L(2n+1,\eta)}$ and its subsequent isotopies as in Figures~\ref{fig:linkinvolution} and \ref{fig:linkquotient} to obtain the corresponding framings $\hat{\lambda}_i$ and $\hat{\mu}_i$ on the $i$th boundary components $\widehat{\bdry N(L_i)}$ of $\bdry (B_{2n+1}, t_L)$.  The picture of $(B_{2n+1}, t_L)$ for general $n$ with framings is shown in Figure~\ref{fig:linktangleisotopy}.

\begin{figure}
\centering
\input{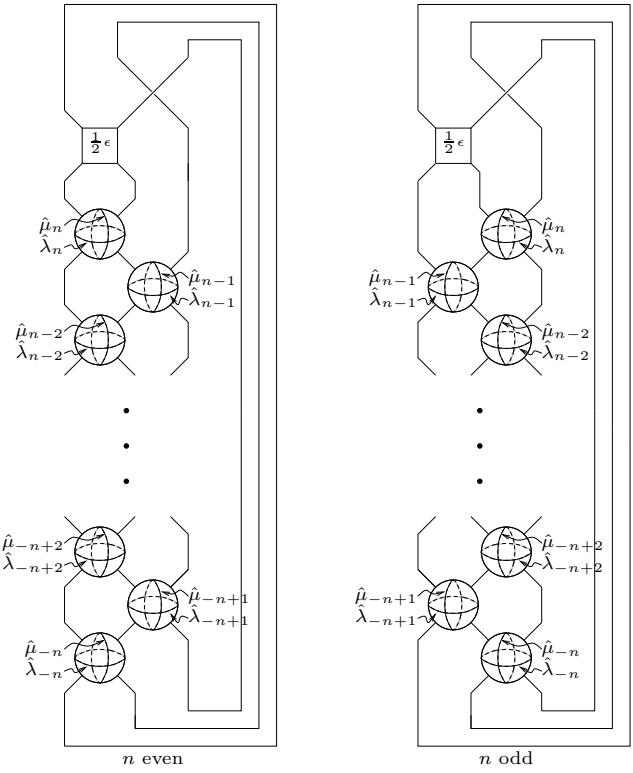}
\caption{The tangles $(B_{2n+1},t_L)$ for each $n$ even and odd.}
\label{fig:linktangleisotopy}
\end{figure}

The minimally twisted chain links $C(2n+1)$ all admit strong involutions too as shown in Figure~\ref{fig:chaininvolution} (a) for $n=2$.  Quotienting the chain link exterior $S^3_{C(2n+1)}$ by the involution yields the tangle $(B_{2n+1},t_C)$ where $B_{2n+1}$ is a $2n+1$ punctured $S^3$ and $t_C$ is the image of the involution axis.  Hence the double branched cover of $B_{2n+1}$ branched over $t_C$ is $S^3_{C(2n+1)}$.  Isotopies of $(B_{2n+1}, t_C)$ for $n=2$ into a ``nice'' form are shown in Figures~\ref{fig:chaininvolution} (b) and (c).  

\begin{figure}
\centering
\input{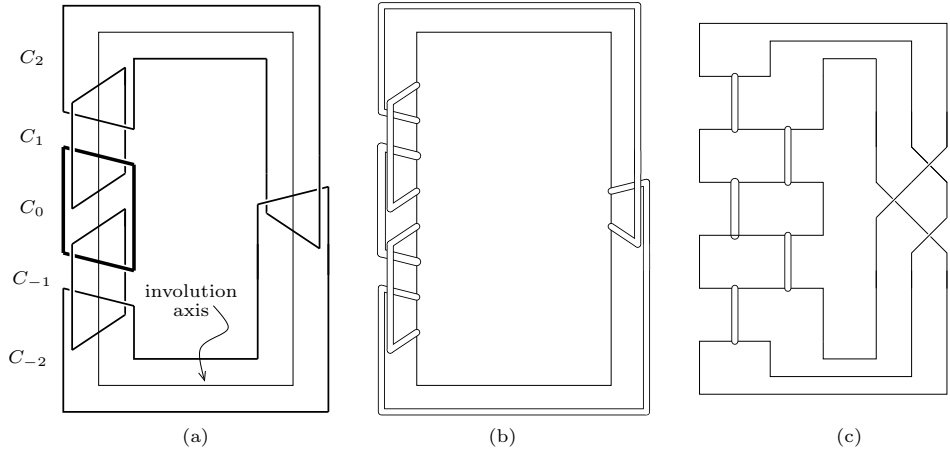}
\caption{The chain link $C(5)$ with an axis of strong involution and the quotient of its complement, $(B_5,t_C)$ followed by an isotopy.}
\label{fig:chaininvolution}
\end{figure}

\comment{
\begin{figure}
\centering
\input{Figures/MT5C-tangleA.pstex_t}
\caption{The chain link $C(5)$ with an axis of strong involution and the quotient of its complement, $(B_5,t_C)$.}
\label{fig:chaininvolution}
\end{figure}

\begin{figure}
\centering
\input{Figures/MT5C-tangleB.pstex_t}
\caption{Isotopies of $(B_5,t_C)$.}
\label{fig:chainquotient}
\end{figure}
}

We may trace the standard longitudes $l_i$ and meridians $m_i$ for $\bdry N(C_i)$ through the quotient of $S^3_{C(2n+1)}$ and the subsequent isotopies as in Figure~\ref{fig:chaininvolution} to obtain the corresponding framings $\hat{l}_i$ and $\hat{m}_i$ on the $i$th boundary components $\widehat{\bdry N(C_i)}$ of $\bdry (B_{2n+1}, t_C)$.  The picture of $(B_{2n+1}, t_C)$ for general $n$ with framings is shown in Figure~\ref{fig:chaintangleisotopy}.

\begin{figure}
\centering
\input{Figures/chaintangleisotopy-framed.pstex_t}
\caption{The tangles $(B_{2n+1}, t_C)$ for each $n$ even and odd.}
\label{fig:chaintangleisotopy}
\end{figure}

\noindent {\em Case 1:} $n$ is even.

Figure~\ref{fig:commontangles} indicates a homeomorphism $\hat{h}$ between the two tangles $(B_{2n+1},t_L)$ and $(B_{2n+1},t_C)$.  The $n$th and $-n$th boundary components of $(B_{2n+1},t_L)$ may absorb the extra twists into their framings as shown in Figure~\ref{fig:eventwisting} to complete the homeomorphism.

\begin{figure}
\centering
\input{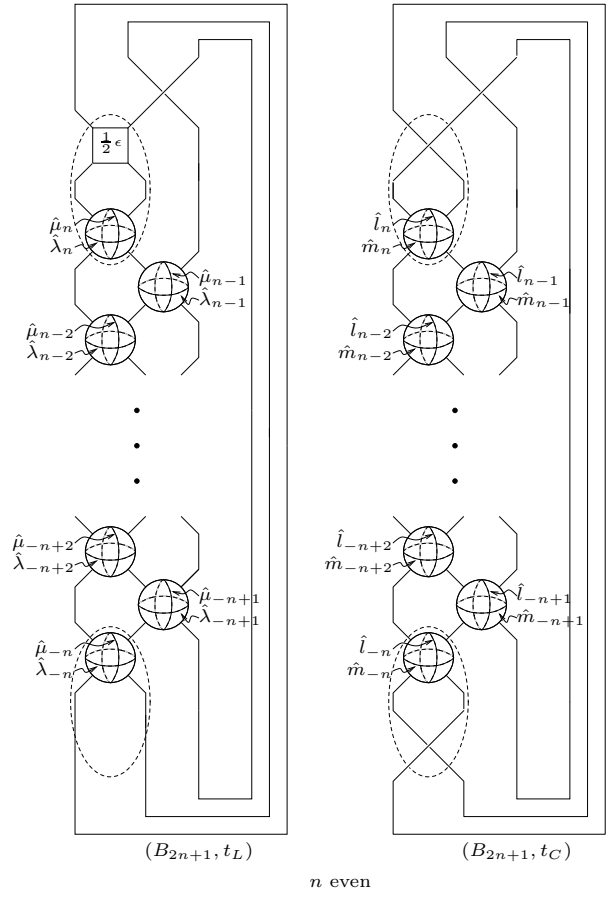}
\caption{The tangles $(B_{2n+1}, t_L)$ and $(B_{2n+1}, t_C)$ for $n$ even.}
\label{fig:commontangles}
\end{figure}

\begin{figure}
\centering
\input{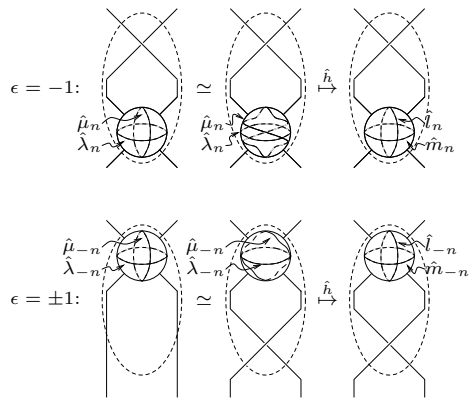}
\caption{Twistings of the $n$th and $-n$th boundary components as needed for $\hat{h}$ for $n$ even.} 
\label{fig:eventwisting}
\end{figure}

The homeomorphism 
\[\hat{h} \colon (B_{2n+1},t_L) \to (B_{2n+1},t_C)\]
then lifts to a homeomorphism
\[h \colon S^3_{L(2n+1,\eta)} \to S^3_{C(2n+1)}\]
of the double branched covers.

\noindent {\em Case 2:} $n$ is odd.

The result of further isotopies of the tangles $(B_{2n+1},t_L)$ and $(B_{2n+1},t_C)$ as depicted in Figures~\ref{fig:linktangleisotopy} and \ref{fig:chaintangleisotopy} for $n$ odd are shown together in Figure~\ref{fig:oddisotopydone} suggesting a homeomorphism $\hat{h}$.  The isotopy of $(B_{2n+1},t_C)$ employs a ``braid relation move.'' 
\begin{figure}
\centering
\input{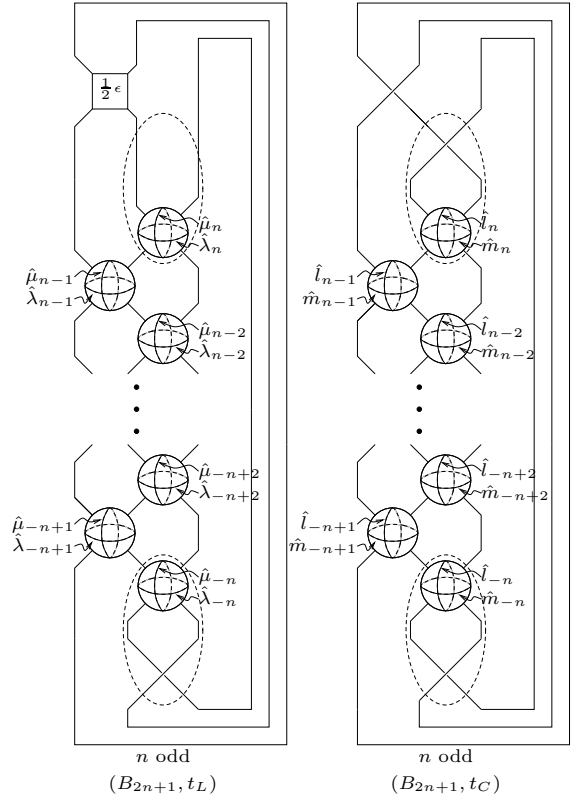}
\caption{The tangles $(B_{2n+1}, t_L)$ and $(B_{2n+1}, t_C)$ for $n$ odd.}
\label{fig:oddisotopydone}
\end{figure}

If $\epsilon = +1$, the $n$th boundary component may absorb an extra twist into its framing as shown in Figure~\ref{fig:oddtwisting+} to complete the homeomorphism $\hat{h}$ between $(B_{2n+1},t_L)$ and $(B_{2n+1},t_C)$ indicated in Figure~\ref{fig:oddisotopydone}.

\begin{figure}
\centering
\input{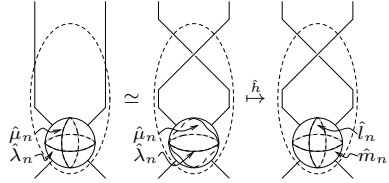}
\caption{Twisting of the $n$th boundary component as needed for $\hat{h}$ for $n$ odd and $\epsilon = +1$.}
\label{fig:oddtwisting+}
\end{figure}

The homeomorphism 
\[\hat{h} \colon (B_{2n+1},t_L) \to (B_{2n+1},t_C)\]
then lifts to a homeomorphism
\[h \colon S^3_{L(2n+1,\eta)} \to S^3_{C(2n+1)}\]
of the double branched covers.

If $\epsilon = -1$, take the reflection of $(B_{2n+1},t_C)$ through the plane of the page.  The resulting tangle $(B_{2n+1}, t_{-C})$ has the complement of $-C(2n+1)$ as its double branched cover.  The $n$th and $-n$th boundary components of $(B_{2n+1},t_L)$ may absorb the extra twists into their framings as shown in Figure~\ref{fig:oddtwisting-} to complete the homeomorphism $\hat{h}$ between $(B_{2n+1},t_L)$ and $(B_{2n+1},t_{-C})$.

\begin{figure}
\centering
\input{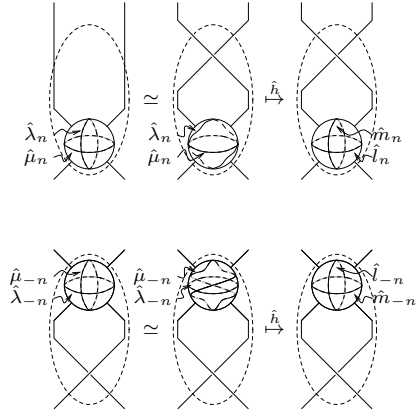}
\caption{Twisting of the $n$th and $-n$th boundary components as needed for $\hat{h}$ for $n$ odd and $\epsilon = -1$.}
\label{fig:oddtwisting-}
\end{figure}

The homeomorphism 
\[\hat{h} \colon (B_{2n+1},t_L) \to (B_{2n+1},t_{-C})\]
then lifts to a homeomorphism
\[h \colon S^3_{L(2n+1,\eta)} \to S^3_{-C(2n+1)}\]
of the double branched covers.

\end{proof}

\subsection{Surgeries on chain links.}

The homeomorphisms $\hat{h}$ in the above proof of Theorem~\ref{thm:homeomorphism} describe how the framed boundary components of $(B_{2n+1},t_L)$ map to the framed boundary components of $(B_{2n+1}, t_{\pm C})$.  In the lift to the double branched covers, this confers how the homeomorphism $h$ maps the boundary components of $S^3_{L(2n+1,\eta)}$ to the boundary components of $S^3_{\pm C(2n+1)}$ in terms of their corresponding meridian-longitude basis pairs.  We thereby obtain mappings of slopes and translate surgery coefficients for $L(2n+1,\eta)$ to surgery coefficients for $\pm C(2n+1)$.

In the following maps, for each $i \in \{-n, \dots, n\}$ the curves $\mu_i$ and $\lambda_i$ on $\bdry N(L_i)$ and the curves $m_i$ and $l_i$ on $\bdry N(C_i)$ are thought of synonymously with their homological representatives so that we may add and subtract them to produce other curves on these tori.  Furthermore these maps should be taken up to reversal of orientations of pairs $\{m_i, l_i\}$.

For each $\hat{h}$ and for each $i \neq \pm n$, the framing curves $\hat{\mu}_i$ and $\hat{\lambda}_i$ map to the framing curves $\hat{l}_i$ and $\hat{m}_i$ respectively.  Lifting to the double branched covers,
\[
\begin{array}{rcl}
\mu_i &\overset{h}{\mapsto}& -l_i \\
\lambda_i &\overset{h}{\mapsto}& m_i.
\end{array}
\]

If $n$ is even, then $\hat{h}$ maps the framings of the $n$th and $-n$th boundary components with twisting as indicated in Figure~\ref{fig:eventwisting}.  In the double branched cover,
\[
\begin{array}{rclcrcl}\label{framing on L_n}
\mu_{-n}+\lambda_{-n} &\overset{h}{\mapsto}&  -l_{-n} &\mbox{ and }&  \mu_n + (1-\epsilon) \lambda_n &\overset{h}{\mapsto} &-l_n \\
\lambda_{-n} & \overset{h}{\mapsto} & m_{-n} && \lambda_n &\overset{h}{\mapsto}& m_n.
\end{array}
\]

If $n$ is odd and $\epsilon = +1$, then $\hat{h}$ maps the framing of the $n$th boundary component with twisting as indicated in Figure~\ref{fig:oddtwisting+}.  The $-n$th boundary component behaves as the others.  In the double branched cover,
\[ 
\begin{array}{rclcrcl}
\mu_{-n} &\overset{h}{\mapsto}& -l_{-n} &\mbox{ and }&  \mu_n - \lambda_n &\overset{h}{\mapsto}& -l_n \\
\lambda_{-n} &\overset{h}{\mapsto}& m_{-n} && \lambda_n &\overset{h}{\mapsto}& m_n.
\end{array}
\] 

If $n$ is odd and $\epsilon = -1$, then $\hat{h}$ maps the framings of the $n$th and $-n$th boundary components with twisting as indicated in Figure~\ref{fig:oddtwisting-}.  In the double branched cover,
\[
\begin{array}{rclcrcl}
\mu_{-n}+2\lambda_{-n} &\overset{h}{\mapsto}&  -l_{-n} &\mbox{ and }& \mu_n + \lambda_n &\overset{h}{\mapsto} &-l_n \\
\lambda_{-n} & \overset{h}{\mapsto} & m_{-n} && \lambda_n &\overset{h}{\mapsto}& m_n.
\end{array}
\]

Therefore, given the slope vector 
\[
\overline{\rho} =\left( (\frac{p}{q})_{-n},(\frac{p}{q})_{-n+1}, \dots,(\frac{p}{q})_{n-1}, (\frac{p}{q})_{n} \right)
\]
for $L(2n+1,\eta)$ in terms of the bases $\{\mu_i, \lambda_i\}$, we may obtain the corresponding slope vector $h(\overline{\rho})$ for $\pm C(2n+1)$ in terms of the bases $\{m_i, l_i\}$ as follows.

If $n$ is even, we obtain
\[
h(\overline{\rho}) =\left( (-\frac{q}{p}+1)_{-n},(-\frac{q}{p})_{-n+1}, \dots,(-\frac{q}{p})_{n-1}, (-\frac{q}{p}+(1-\epsilon))_{n} \right)
\]
for $C(2n+1)$.

If $n$ is odd, we obtain
\[
h(\overline{\rho}) =\left( (-\frac{q}{p}+(1-\epsilon))_{-n},(-\frac{q}{p})_{-n+1}, \dots,(-\frac{q}{p})_{n-1}, (-\frac{q}{p}-\epsilon)_{n} \right)
\]
for $\epsilon \cdot C(2n+1)$.

In particular, given $K \in \K_\eta$ with complement $S^3_K$ described as $\overline{\rho}$-Dehn surgery on $L(2n+1,\eta)$ where
\[
\overline{\rho} = \left( (-\frac{1}{r_n})_{-n}, (-\frac{1}{r_{n-1}})_{-n+1}, \dots, (-\frac{1}{r_1})_{-1}, (\emptyset)_{0}, (\frac{1}{r_1})_{1}, \dots, (\frac{1}{r_{n-1}})_{n-1}, (\frac{1}{r_n})_{n} \right)
\]
as in Proposition~\ref{prop:K_as_surgery} (so that the coefficient $\emptyset$ indicates that its component is to be left unfilled), then for $n$ even
\begin{multline*}
h(\overline{\rho}) = \left( (r_n +1)_{-n}, (r_{n-1})_{-n+1}, \dots, (r_1)_{-1}, (\emptyset)_{0},\right.\\ 
\left. (-r_1)_{1}, \dots, (-r_{n-1})_{n-1}, (-r_n + (1-\epsilon))_{n} \right)
\end{multline*}
and for $n$ odd
\begin{multline*}
h(\overline{\rho}) = \left( (r_n + (1-\epsilon))_{-n}, (r_{n-1})_{-n+1}, \dots, (r_1)_{-1}, (\emptyset)_{0},\right.\\ 
\left. (-r_1)_{1}, \dots, (-r_{n-1})_{n-1}, (-r_n -\epsilon)_{n} \right).
\end{multline*}
Furthermore, the $0$--slope on the $0$th component of $S^3_{\pm C(2n+1)}(h(\overline{\rho})) \cong S^3_K$ is the meridional slope for $K$.

\begin{figure}
\centering
\input{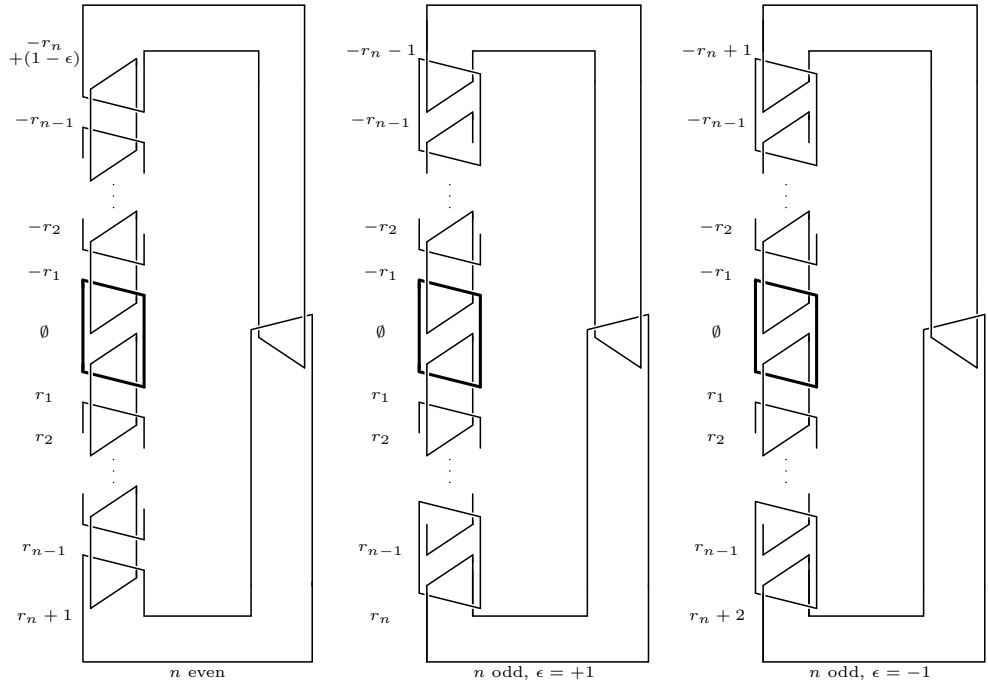}
\caption{Surgery description of $K \in \K$ on $C(2n+1)$.}
\label{fig:MTCsurgery}
\end{figure}

Thus if $K \in \K_\eta$ has slope $[0, r_n, \dots, r_2, r_1]$ with $n$ even or $[r_n, \dots, r_2, r_1]$ with $n$ odd, then $S^3_K$ has a surgery description on $\pm C(2n+1)$ as shown in Figure~\ref{fig:MTCsurgery}.  Again, the $0$--slope on the $0$th component is the meridional slope for $K$ giving the $S^3$ surgery.

This proves 
\begin{thm}\label{thm:chainlinksurgdesc}
The knots in Berge's families (VII) and (VIII) all have surgery descriptions on some minimally twisted chain link, $\pm C(2n+1)$ for some $n \in \N$.
\end{thm}

\subsection{Hyperbolicity of the links $L(2n+1, \eta)$.}

Neumann-Reid give the following theorem from which we derive the subsequent corollary.
\begin{thm}[Theorem 5.1 (ii) of \cite{neumannreid}]\label{thm:neumannreid}
$C(p,s)$ has a hyperbolic structure (complete of finite volume) if and only if $\{|p+s|,|s|\} \not\subset \{0,1,2\}$.
\end{thm}

\begin{cor}\label{cor:chainlinkhyp}
The links $C(2n+1)$ are hyperbolic for $n > 1$.
\end{cor}

\begin{proof}
Our links $C(2n+1)$ correspond to the links ${C(2n+1,-n-1)}$ or ${C(2n+1, -n)}$ of \cite{neumannreid} depending on how the last clasping is done.  Since $n > 1$ Theorem~\ref{thm:neumannreid} applies.
\end{proof}

With the homeomorphism of link exteriors of Theorem~\ref{thm:homeomorphism} and Corollary~\ref{cor:chainlinkhyp}, we conclude
\begin{cor}\label{cor:Lishyp}
$L(2n+1,\eta)$ is a hyperbolic link for $n\geq2$.
\end{cor}
\bibliography{MathBiblio}
\bibliographystyle{plain}
\end{document}